\documentclass{article}
\usepackage{amssymb}
\usepackage{amsmath}
\usepackage{latexsym}
\usepackage{longtable}

\input{prepictex}
\input{pictex}
\input{postpictex}
\begin{document}
\bibliographystyle{plain}
\baselineskip=14pt
\newtheorem{lemma}{Lemma}[section]
\newtheorem{theorem}[lemma]{Theorem}
\newtheorem{prop}[lemma]{Proposition}
\newtheorem{cor}[lemma]{Corollary}
\newtheorem{definition}[lemma]{Definition}
\newtheorem{definitions}[lemma]{Definitions}
\newtheorem{remark}[lemma]{Remark}
\newtheorem{conj}[lemma]{Conjecture}
\renewcommand{\arraystretch}{2}
\newcommand{\TF}{\widetilde{F}}
\newcommand{\C}{{\mathbb{C}}}
\newcommand{\N}{{\mathbb{N}}}
\newcommand{\Q}{{\mathbb{Q}}}
\newcommand{\R}{{\mathbb{R}}}
\newcommand{\Z}{{\mathbb{Z}}}
\newcommand{\B}{{\mathbf{B}}}
\newcommand{\CA}{{\mathcal{A}}}
\newcommand{\CL}{{\mathcal{L}}}
\newcommand{\0}{-\,}
\newcommand{\cdi}{\mbox{CD$(\mathbf{i})$}}
\newcommand{\cdii}{\mbox{CD$(\mathbf{i}')$}}
\begin{center}
\baselineskip=11pt
{\LARGE Canonical Bases and Piecewise-linear Combinatorics}

\vspace*{1cm}

{\Large Roger Carter} \\
\vspace*{0.5cm}
{\em Mathematics Institute, University of Warwick, Coventry CV4 7AL, England}
\\
{E-mail: rwc@maths.warwick.ac.uk}

\vspace*{0.3cm}
{\Large Robert Marsh} \\
\vspace*{0.5cm}
{\em Department of Mathematics and Computer Science, University of Leicester,
University Road, Leicester LE1 7RH, England} \\
{E-mail: R.Marsh@mcs.le.ac.uk}

\baselineskip=22pt

{\bf Abstract} \\
\parbox[t]{5.15in}{
Let $U_q$ be the quantum group associated to a Lie algebra
$\bf g$ of rank $n$. The negative part $U_q^-$ of $U_q$ has
a canonical basis $\mathbf{B}$ with favourable properties
(see Kashiwara~\cite{kash2} and Lusztig~\cite[\S14.4.6]{lusztig6}). The
approaches of Lusztig and Kashiwara lead to a set of
alternative parametrizations of the canonical basis, one for each reduced
expression for the longest word in the Weyl group of $\bf g$.
We describe the authors' recent work establishing close relationships
between the Lusztig cones, canonical basis elements and the regions of
linearity of reparametrization functions arising from the above
parametrizations in type $A_4$ and give some speculations for type $A_n$.

} \\ \ \\
\end{center}

\section{Introduction}

Let $\mathbf{g}$ be a finite-dimensional simple Lie algebra over ${\mathbb{C}}$
and $U_q(\mathbf{g})$ be the corresponding quantized enveloping algebra over
${\mathbb{C}}[v,v^{-1}]$. Let
$$U_q(\mathbf{g})=U_q^- \otimes U_q^0 \otimes U_q^+$$
be the triangular decomposition of $U_q(\mathbf{g})$. Let $\mathbf{B}$ be the
canonical basis of $U_q^-$, introduced independently by Lusztig and Kashiwara. It
is natural to ask how the elements of $\mathbf{B}$ are expressed in terms of the
generators $F_1,F_2,\ldots ,F_n$ of $U_q^-$. This is a difficult question
which is known completely in only a few low rank cases, but the attempt to
understand it has led to a remarkable theory of piecewise-linear combinatorics
associated with the canonical basis.

In this paper we consider only the case in which $\mathbf{g}$ has type $A_n$.
We describe the situation for small values of $n$ before giving some
speculations for arbitrary $n$.

First suppose that $\mathbf{g}$ has type $A_1$. We write:
\begin{eqnarray*}
{\mathbb{N}} & = & \{a\in {\mathbb{Z}}\,:\,a\geq 0\}, \\
\left[ a \right]  & = & \frac{v^a-v^{-a}}{v-v^{-1}},\ \ \mbox{for\ }a\in {\mathbb{N}},  \\
\left[ a \right]! & = & [1][2]\cdots [a],\ \ \mbox{for\ }a\in {\mathbb{N}},
\end{eqnarray*}
and $F_1^{(a)}=F_1^a/[a]!$.

It was shown by Lusztig~\cite{lusztig2} that $\mathbf{B}=\{F_1^{(a)}\,:\,
a\in {\mathbb{N}}\}$. Thus the canonical basis elements are the
quantized divided powers of the generator $F_1$.

Next, suppose that $\mathbf{g}$ has type $A_2$. The canonical basis in this case
was also determined by Lusztig~\cite{lusztig2}. We have
$$\begin{array}{rccl}
\mathbf{B} & = & & \{F_1^{(a)}F_2^{(b)}F_1^{(c)}\,:\,\ a,b,c\in {\mathbb{N}},\ 
b\geq a+c\} \\
        &   & & \cup\{F_2^{(a)}F_1^{(b)}F_2^{(c)}\,:\,\ a,b,c\in {\mathbb{N}},\ 
b\geq a+c\}.
\end{array}$$
When $b=a+c$ one has the relation
$$F_1^{(a)}F_2^{(a+c)}F_1^{(c)}=F_2^{(c)}F_1^{(a+c)}F_2^{(a)},$$
but apart from this the above elements of $\mathbf{B}$ are all distinct. Thus in
this case each canonical basis element may be written as a monomial in the
generators $F_1,F_2$. The two types of monomials which arise are related to
the two reduced decompositions of the longest element $w_0$ of the Weyl group
$W$ of $\mathbf{g}$. We have $W=\langle s_1,s_2 \rangle$ and
$w_0=s_1s_2s_1=s_2s_1s_2$. Each of these reduced words for $w_0$ gives rise
to a type of monomial in the canonical basis.

We now turn to the case when $\mathbf{g}$ has type $A_3$. Lusztig~\cite{lusztig9}
obtained many, but not all, elements of $\mathbf{B}$ as monomials in the
generators $F_1,F_2,F_3$ and gave an example of an element of $\mathbf{B}$ which
could not be written as a monomial in $F_1,F_2,F_3$. The remaining elements
of $\mathbf{B}$ were determined by Xi~\cite{xi3}. In order to describe the
monomials in $\mathbf{B}$ we consider reduced words for $w_0$. In type $A_3$
$w_0$ has $16$ reduced expressions. However it is natural to divide them
into equivalence classes called commutation classes; two reduced words being
equivalent if one can be obtained from the other by a succession of
commutations in the Coxeter generators of $W$. For example,
$s_1s_3s_2s_1s_3s_2$ is in the same commutation class as $s_3s_1s_2s_3s_1s_2$,
where the generators are numbered as in the Dynkin diagram (see Figure $1$).

\begin{figure}[htbp]

\beginpicture

\setcoordinatesystem units <1cm,1cm>             
\setplotarea x from -4 to 6, y from 0 to 3       


\put {Figure $1$: Dynkin diagram of type $A_3$.}[c] at 4 0.5

\scriptsize{

\multiput {$\circ$} at 3   2 *2 1 0 /      %

\linethickness=1pt           

\putrule from 3.05 2 to 3.95 2  %
\putrule from 4.05 2 to 4.95 2  

\put {$1$}   [l] at 2.9 1.75
\put {$2$}   [l] at 3.9 1.75
\put {$3$}   [l] at 4.9 1.75

}

\endpicture

\end{figure}

There are $8$ commutation classes of reduced words for $w_0$ in type $A_3$.
Each of these gives rise to a family of monomials in $F_1,F_2,F_3$ which
lie in $\mathbf{B}$. For example, the reduced word $w_0=s_1s_3s_2s_1s_3s_2$ gives
rise to monomials
$$F_1^{(a)}F_3^{(b)}F_2^{(c)}F_1^{(d)}F_3^{(e)}F_2^{(f)},\ \ \ \ 
a,b,c,d,e,f\in {\mathbb{N}},$$
and such a monomial was shown by Lusztig to lie in the canonical basis $\mathbf{B}$
provided
$$c\geq a+d,\ \ c\geq b+e,\ \ d+e\geq c+f.$$
These inequalities come from considering two consecutive occurrences of a
given Coxeter generator $s_i$ in the given reduced word. The sum of the
exponents corresponding to these occurrences of $s_i$ is less than or equal to
the sum of the exponents corresponding to Coxeter generators $s_j$ between
these two occurrences of $s_i$ such that $s_j$ does not commute with $s_i$.
Lusztig obtained $8$ families of monomials in $\mathbf{B}$ in this way,
and the remaining elements of $\mathbf{B}$ determined by Xi are linear
combinations of monomials with coefficients which are quantum binomial
coefficients.

We now suppose that $\mathbf{g}$ has type $A_4$. Here the situation is more
involved and has been investigated by the authors. We shall outline the
situation in the present paper, and hope to publish the proofs in a
subsequent article. For each reduced word $w_0=s_{i_1}s_{i_2}\cdots s_{i_k},$
where $k=\ell (w_0)$, we write $\mathbf{i}=(i_1,i_2,\ldots ,i_k)$. For each
such $\mathbf{i}$, we define a subset $C_\mathbf{i}$ of ${\mathbb{N}}^k$ whose definition
is motivated by the rule mentioned above for a monomial to lie in $\mathbf{B}$ in
type $A_3$. We define $C_\mathbf{i}$ to be the set of those
$\mathbf{a}\in {\mathbb{N}}^k$ such that for each pair $t,t'\in [1,k]$ with $t<t'$,
$i_t=i_{t'}$, $i_p\not=i_t$ for $t<p<t'$, we have
$\sum_p a_p\geq a_t+a_{t'}$, summed over all $p$ with $t<p<t'$ such that $s_p$
does not commute with $s_{i_t}$.

The cone $C_\mathbf{i}$ will be called the {\em Lusztig cone}
associated with $\mathbf{i}$.
The elements of $C_\mathbf{i}$ give rise to monomials. Let
$$M_\mathbf{i}=\{F_{i_1}^{(a_1)}F_{i_2}^{(a_2)}\cdots F_{i_k}^{(a_k)}\,:\,
\mathbf{a}=(a_1,a_2,\ldots ,a_k)\in C_\mathbf{i}\}.$$
It was shown by Marsh~\cite{me7} that when $\mathbf{g}$ has type $A_4$ we have
$M_\mathbf{i}\subseteq \mathbf{B}$ for each $\mathbf{i}$. In type $A_4$ there are $62$
commutation classes of reduced words for $w_0$, and we obtain in this way
$62$ families of monomials in the canonical basis. These are far from being
the only elements of $\mathbf{B}$, however, and the remaining elements are not
known as expressions in terms of $F_1,F_2,F_3,F_4$.

\section{PBW-type Bases and the Canonical Basis}

We shall now recall Lusztig's approach to the canonical basis. For each
reduced word $\mathbf{i}$ for $w_0$, Lusztig defined~\cite{lusztig2} a PBW-type
basis $B_\mathbf{i}$ of $U_q^-$. The reduced word $w_0=s_{i_1}s_{i_2}\cdots
s_{i_k}$ gives rise to a total order on the set $\Phi^+$ of positive roots
of $\mathbf{g}$. We have
$$\Phi^+=\{\alpha^1,\alpha^2,\ldots ,\alpha^k\},$$
where
$$\alpha^1=\alpha_{i_1},\ \alpha^2=s_{i_1}(\alpha_{i_2}),\ldots \alpha^k=
s_{i_1}s_{i_2}\cdots s_{i_{k-1}}(\alpha_{i_k}),$$
and $\alpha_1,\alpha_2,\ldots ,\alpha_n$ are the simple roots. By using
braid group operations, Lusztig defines for each $\alpha\in \Phi^+$ a root
vector $F_{\alpha}\in U_q^-$, beginning with $F_{\alpha_i}=F_i$. We write
$F_\mathbf{i}^\mathbf{c}=F_{\alpha^1}^{(c_1)}F_{\alpha^2}^{(c_2)}\cdots
F_{\alpha^k}^{(c_k)}$ for $\mathbf{c}=(c_1,c_2,\ldots ,c_k)\in {\mathbb{N}}^k$.
Then the set $B_\mathbf{i}=\{F_\mathbf{i}^\mathbf{c}\,:\,\mathbf{c}\in {\mathbb{N}}^k\}$ is a
basis for $U_q^-$, called the PBW-type basis associated to $\mathbf{i}$.

The lattice $\mathcal{L}={\mathbb{C}}[v]B_\mathbf{i}$ was shown by Lusztig to be
independent of $\mathbf{i}$ and there is a bijective map
\begin{eqnarray*}
\mathbf{B} & \longrightarrow & B_\mathbf{i} \\
b & \mapsto & F_\mathbf{i}^\mathbf{c}
\end{eqnarray*}
such that $b\equiv F_\mathbf{i}^\mathbf{c}\mod v \mathcal{L}.$
We write $\varphi_\mathbf{i}(b)=\mathbf{c}$. Then the map $\varphi_{\mathbf{i}}:
\mathbf{B}\rightarrow {\mathbb{N}}^k$ is bijective. This gives, for each
$\mathbf{i}$, a parametrization of $\mathbf{B}$ by elements of ${\mathbb{N}}^k$.

Lusztig introduced~\cite{lusztig6} two particular reduced words $\mathbf{j}$,
$\mathbf{j'}$ for $w_0$. In type $A_n$ these are as follows. Suppose the
vertices of the Dynkin diagram of $A_n$ are labelled as in Figure $2$.

\begin{figure}[htbp]

\beginpicture

\setcoordinatesystem units <1cm,1cm>             
\setplotarea x from -2.5 to 12, y from 0 to 3       

\put {Figure $2$: Dynkin diagram of type $A_n$.}[c] at 5.5 0.5

\scriptsize{

\multiput {$\circ$} at 3   2 *1 1 0 /      %
\multiput {$\circ$} at 7   2 *1 1 0 /      

\linethickness=1pt           

\putrule from 3.05 2 to 3.95 2  %
\putrule from 7.05 2 to 7.95 2  

\setdashes <2mm,1mm>          %
\putrule from 4.05 2 to 6.95 2  

\put {$1$}   [l] at 2.9 1.75
\put {$2$}   [l] at 3.9 1.75
\put {$n-1$} [l] at 6.65 1.75
\put {$n$} [l] at 7.9 1.75

}

\endpicture

\end{figure}

Let $\mathbf{j}$ be $135\cdots 246\cdots 135\cdots 246\cdots $, where
$k=\frac{1}{2}n(n+1)$ factors are taken, and let $\mathbf{j'}$ be \linebreak
$246\cdots 135\cdots 246\cdots 135\cdots $, where again $k$ factors are taken.
Then $\mathbf{j},\mathbf{j'}$ are reduced words for $w_0$ and their
commutation classes are as far apart as possible, in the sense that they
give opposite orderings on $\Phi^+$. Consider the parametrizations
$\varphi_{\mathbf{j}}:\mathbf{B}\rightarrow {\mathbb{N}}^k$ and
$\varphi_{\mathbf{j'}}:\mathbf{B}\rightarrow {\mathbb{N}}^k$ of $\mathbf{B}$
corresponding to $\mathbf{j},\mathbf{j'}$ and let
$R=\varphi_{\mathbf{j'}}\varphi_{\mathbf{j}}^{-1}:
{\mathbb{N}}^k\rightarrow {\mathbb{N}}^k$ be the bijective map which relates
them.

The function $R$ is the restriction to ${\mathbb{N}}^k$ of a piecewise-linear
map $R:{\mathbb{R}}^k\rightarrow {\mathbb{R}}^k$ and the results in low rank
cases suggest that the regions of linearity of $R$ are related to the
different types of canonical basis elements.

For example, in type $A_2$, we have $\mathbf{j}=121$, $\mathbf{j'}=212$,
and $R$ has two regions of linearity. These correspond to the two types of
monomial in the canonical basis.

In type $A_3$, the function $R$ has $10$ regions of linearity, and $8$ of
these regions correspond to the $8$ families of monomials in the canonical
basis. The remaining $2$ regions correspond to the remaining non-monomial
elements of $\mathbf{B}$ obtained by Xi. These two regions may be
distinguished from the other $8$ regions as follows.

In type $A_3$, we have $R:{\mathbb{R}}^6\rightarrow {\mathbb{R}}^6$. The $8$ regions
of linearity of $R$ which give rise to monomials in $\mathbf{B}$ are each defined
by $3$ inequalities whereas the $2$ remaining regions of linearity are
each defined by $4$ inequalities. This suggests that regions of linearity
defined by the minimum number of inequalities might give rise to canonical
basis elements of a particularly favourable form.

This turns out to be the case in type $A_4$ also. This time we have a
piecewise-linear map $R:{\mathbb{R}}^{10}\rightarrow {\mathbb{R}}^{10}$, which was shown
by Carter to have 144 regions of linearity. Of these regions, $62$ are
defined by $6$ inequalities, $70$ by $7$ inequalities, $10$ by $8$
inequalities, and $2$ by $11$ inequalities. It is striking that the number
of regions defined by the minimum number of inequalities is equal to the
number of commutation classes of reduced words for $w_0$.

The authors have shown that in type $A_4$ there is a natural bijection between
commutation classes of reduced words for $w_0$ and regions of linearity for
$R$ defined by the minimum number of inequalities. This comes about as
follows. For each reduced word $\mathbf{i}$ for $w_0$ we have a corresponding
family $M_{\mathbf{i}}$ of monomials in $\mathbf{B}$, as described above.
We consider the set of points in ${\mathbb{N}}^{10}$ which parametrize these
elements of $\mathbf{B}$ under the map $\varphi_{\mathbf{j}}$: let
$X_{\mathbf{i}}^+=\varphi_{\mathbf{j}}(M_{\mathbf{i}})$. Then one can show: \\
(a) $X_{\mathbf{i}}^+$ is the set of all points in a region of linearity
$X_{\mathbf{i}}$ of $R$ with coordinates in $\mathbb{N}$. \\
(b) $X_{\mathbf{i}}$ is a region defined by the minimum number of inequalities. \\
(c) The map $\mathbf{i}\rightarrow X_{\mathbf{i}}$ is a bijection between
commutation classes of reduced words for $w_0$ and regions of $R$ defined
by the minimum number of inequalities.

In fact, both the set of commutation classes of reduced words for $w_0$ and
the set of regions of $R$ defined by the minimum number of inequalities can
be given a natural graph structure, and the map $\mathbf{i}\rightarrow
X_{\mathbf{i}}$ is then an isomorphism of graphs.

\section{Rectangle Calculus}
It is natural to ask whether the region $X_{\mathbf{i}}$ can be described in
terms of $\mathbf{i}$ without using ideas concerned with the canonical basis,
but simply in combinatorial terms. This can be done by a form of
combinatorics which we call the rectangle calculus. The basic idea is to
associate with the reduced word $\mathbf{i}$ a set of linearly independent
vectors such that $X_{\mathbf{i}}$ is the set of all linear combinations of
these vectors with non-negative coefficients. These vectors will be called
{\em spanning vectors} of $X_{\mathbf{i}}$.

We first introduce the idea of a partial quiver. This is a Dynkin diagram with
arrows on certain edges, such that the set of edges with arrows in non-empty
and connected. An example of a partial quiver of type $A_7$ is given in Figure
$3$.

\begin{figure}[htbp]

\beginpicture

\setcoordinatesystem units <1cm,1cm>             
\setplotarea x from -2 to 12, y from 0 to 3       

\put {Figure $3$: A partial quiver of type $A_7$.}[c] at 6 0.5

\scriptsize{

\multiput {$\circ$} at 3   2 *6 1 0 /      %

\putrule from 3.05 2 to 3.95 2  %
\putrule from 4.05 2 to 4.95 2  %
\putrule from 5.05 2 to 5.95 2  %
\putrule from 6.05 2 to 6.95 2  %
\putrule from 7.05 2 to 7.95 2  %
\putrule from 8.05 2 to 8.95 2  %

\setlinear \plot  5.4 2.1  5.5 2  / %
\setlinear \plot  5.4 1.9  5.5 2  / %
\setlinear \plot  6.6 2.1  6.5 2  / %
\setlinear \plot  6.6 1.9  6.5 2  / %
\setlinear \plot  7.6 2.1  7.5 2  / %
\setlinear \plot  7.6 1.9  7.5 2  / %

}

\endpicture

\end{figure}

We are concerned with partial quivers of type $A_n$, and shall use $L$ or
$R$ to indicate whether an arrow goes left or right. Thus the above partial
quiver is denoted $\0\0RLL-$.

The possible partial quivers of type $A_4$ are $LLL$, $RLL$, $LRL$,
$LLR$, $LRR$, $RLR$, $RRL$, $RRR$, $LL-$, $\0LL$, $LR-$, $\0LR$, $RL-$,
$\0RL$, $RR-$, $\0RR$, $L\0-$, $\0L-$, $\0\0L$, $R\0-$, $\0R-$, $\0\0R$.

We shall now explain a procedure by which each reduced word $\mathbf{i}$ for
$w_0$ in type $A_n$ gives a set of $\frac{1}{2}n(n-1)$ partial quivers.
We first write down the braid diagram of $\mathbf{i}$. This determines a set
of chambers, and for each bounded chamber we write down the corresponding
chamber set, which is the subset of $\{1,2,\ldots ,n+1\}$ corresponding to
the strings which pass below the chamber. We illustrate this by the example
in which $n=4$ and $w_0=s_2s_3s_4s_3s_1s_2s_1s_3s_2s_4$ --- see Figure $4$.


\begin{figure}[htbp]

\beginpicture

\setcoordinatesystem units <0.8cm,0.6cm>             
\setplotarea x from -4 to 14, y from -4 to 5.5       

\linethickness=0.5pt           

\put{Figure $4$: A chamber diagram.}[c] at 6 -3

\put{$1$}[c] at -0.3 4 %
\put{$2$}[c] at -0.3 3 %
\put{$3$}[c] at -0.3 2 %
\put{$4$}[c] at -0.3 1 %
\put{$5$}[c] at -0.3 0 %

\put{$2$}[c] at 1.5 -1 %
\put{$3$}[c] at 2.5 -1 %
\put{$4$}[c] at 3.5 -1 %
\put{$3$}[c] at 4.5 -1 %
\put{$1$}[c] at 5.5 -1 %
\put{$2$}[c] at 6.5 -1 %
\put{$1$}[c] at 7.5 -1 %
\put{$3$}[c] at 8.5 -1 %
\put{$2$}[c] at 9.5 -1 %
\put{$4$}[c] at 10.5 -1 %

\setlinear \plot 0 0  3 0 / %
\setlinear \plot 3 0  5 2 / %
\setlinear \plot 5 2  6 2 / %
\setlinear \plot 6 2  8 4 / %
\setlinear \plot 8 4  12 4 / %

\setlinear \plot 0 1  2 1 / %
\setlinear \plot 2 1  3 2 / %
\setlinear \plot 3 2  4 2 / %
\setlinear \plot 4 2  5 1 / %
\setlinear \plot 5 1  8 1 / %
\setlinear \plot 8 1  10 3 / %
\setlinear \plot 10 3  12 3 / %

\setlinear \plot 0 2  1 2 / %
\setlinear \plot 1 2  2 3 / %
\setlinear \plot 2 3  5 3 / %
\setlinear \plot 5 3  6 4 / %
\setlinear \plot 6 4  7 4 / %
\setlinear \plot 7 4  8 3 / %
\setlinear \plot 8 3  9 3 / %
\setlinear \plot 9 3  10 2 / %
\setlinear \plot 10 2 12 2 / %

\setlinear \plot 0 3  1 3 / %
\setlinear \plot 1 3  4 0 / %
\setlinear \plot 4 0  10 0 / %
\setlinear \plot 10 0  11 1 / %
\setlinear \plot 11 1  12 1 / %

\setlinear \plot 0 4  5 4 / %
\setlinear \plot 5 4  7 2 / %
\setlinear \plot 7 2  8 2 / %
\setlinear \plot 8 2  9 1 / %
\setlinear \plot 9 1  10 1 / %
\setlinear \plot 10 1 11 0 / %
\setlinear \plot 11 0 12 0 / %

\put{$25$}[c] at 3.5 1.5 %
\put{$245$}[c] at 4 2.5 %
\put{$2$}[c] at 7 0.5 %
\put{$24$}[c] at 6.5 1.5 %
\put{$1245$}[c] at 6.5 3.5 %
\put{$124$}[c] at 8 2.5 %

\endpicture

\end{figure}

Each chamber set
obtained in this way is a subset of $\{1,2,\ldots ,n+1\}$ which is not
an initial or terminal subset, i.e. not of form $\{1,2,\ldots ,i\}$ or
$\{i,i+1,\ldots ,n+1\}$ for any $i$.

Let $\mathcal{S}$ be the set of all subsets of $\{1,2,\ldots ,n+1\}$ which are
not initial or terminal subsets and $\mathcal{P}$ be the set of all
partial quivers of type $A_n$. We shall describe a bijection from
$\mathcal{S}$ to $\mathcal{P}$. We first number the edges of the Dynkin
diagram as shown in Figure $5$.

\begin{figure}[htbp]

\beginpicture

\setcoordinatesystem units <1cm,1cm>             
\setplotarea x from -2.7 to 12, y from 0 to 3       
\put {Figure $5$: Edge numbering of Dynkin diagram.}[c] at 5.5 0.5

\scriptsize{

\multiput {$\circ$} at 2   2 *2 1 0 /      %
\multiput {$\circ$} at 7   2 *2 1 0 /      

\linethickness=1pt           

\putrule from 2.05 2 to 2.95 2  %
\putrule from 3.05 2 to 3.95 2  %
\putrule from 7.05 2 to 7.95 2  
\putrule from 8.05 2 to 8.95 2  %

\setdashes <2mm,1mm>          %
\putrule from 4.05 2 to 6.95 2  

\put {$1$}   [c] at 9.5 2.2
\put {$2$}   [c] at 8.5 2.2
\put {$3$}   [c] at 7.5 2.2
\put {$n-1$} [c] at 3.5 2.2
\put {$n$}   [c] at 2.5 2.2
\put {$n+1$} [c] at 1.5 2.2

}

\endpicture

\end{figure}

Thus the edges are numbered $2,3,\ldots ,n$ from right to left and $1,n+1$
are regarded as virtual edges.

Let $S\in \mathcal{S}$ and $P$ be the corresponding partial quiver. $P$
is obtained from $S$ by the following rules.

(a) If $1,n+1\not\in S$, the entries in $S$ give the edges of type $L$
in $P$, edges intermediate between those of type $L$ having type $R$. The
leftmost and rightmost labelled edges of $P$ have type $L$. \\
(b) Now suppose $1\in S$. Let $i$ be such that $1,2,\ldots ,i\in S$
but $i+1\not\in S$. Then edge $i+1$ is labelled $R$ and is the rightmost
labelled edge in $P$. \\
(c) Now suppose $n+1\in S$. Let $i$ be such that $i,i+1,\ldots ,n+1$
lie in $S$ but $i-1\not\in S$. Then edge $i-1$ of $P$ is labelled $R$ and
is the leftmost labelled edge in $P$. \\
(d) The elements of $S$ not in an initial segment as in (b) or a
terminal segment as in (c) give rise to edges $L$ or $R$ in $P$ as in (a).

{\bf Example}. \\
Let $n=13$. If $S=\{1,2,3,4,7,8,11\}$, then $P=\0\0LRRLLRR\0-\hspace*{1.3pt}-$. See Figure $6$ for an explanatory diagram.

\begin{figure}[htbp]

\beginpicture

\setcoordinatesystem units <1cm,1cm>             
\setplotarea x from -0.4 to 15.5, y from 0 to 3       
\put {Figure $6$: How to calculate the subset corresponding to a partial quiver.}[c] at 8 0.5

\scriptsize{

\multiput {$\circ$} at 2   2 *12 1 0 /      %

\linethickness=1pt           

\putrule from 2.05 2 to 2.95 2  %
\putrule from 3.05 2 to 3.95 2  %
\putrule from 4.05 2 to 4.95 2  %
\putrule from 5.05 2 to 5.95 2  %
\putrule from 6.05 2 to 6.95 2  %
\putrule from 7.05 2 to 7.95 2  %
\putrule from 8.05 2 to 8.95 2  %
\putrule from 9.05 2 to 9.95 2  %
\putrule from 10.05 2 to 10.95 2  %
\putrule from 11.05 2 to 11.95 2  %
\putrule from 12.05 2 to 12.95 2  %
\putrule from 13.05 2 to 13.95 2  %

\put {$1$}   [c] at 14.5 2.2
\put {$2$}   [c] at 13.5 2.2
\put {$3$}   [c] at 12.5 2.2
\put {$4$}   [c] at 11.5 2.2
\put {$5$}   [c] at 10.5 2.2
\put {$6$}   [c] at 9.5 2.2
\put {$7$}   [c] at 8.5 2.2
\put {$8$}   [c] at 7.5 2.2
\put {$9$}   [c] at 6.5 2.2
\put {$10$}  [c] at 5.5 2.2
\put {$11$}  [c] at 4.5 2.2
\put {$12$}  [c] at 3.5 2.2
\put {$13$}  [c] at 2.5 2.2
\put {$14$}  [c] at 1.5 2.2

\put {$-$} [c] at 2.5 1.5
\put {$-$} [c] at 3.5 1.5
\put {$L$} [c] at 4.5 1.5
\put {$R$} [c] at 5.5 1.5
\put {$R$} [c] at 6.5 1.5
\put {$L$} [c] at 7.5 1.5
\put {$L$} [c] at 8.5 1.5
\put {$R$} [c] at 9.5 1.5
\put {$R$} [c] at 10.5 1.5
\put {$-$} [c] at 11.5 1.5
\put {$-$} [c] at 12.5 1.5
\put {$-$} [c] at 13.5 1.5

}

\endpicture

\end{figure}

Applying this bijection $\mathcal{S}\rightarrow \mathcal{P}$ to the chamber
sets obtained from a reduced word $\mathbf{i}$ we obtain a set of
$\frac{1}{2}n(n-1)$ partial quivers associated with $\mathbf{i}$. In the
above example with $n=4$ and $w_0=s_2s_3s_4s_3s_1s_2s_1s_3s_2s_4$ we obtain
the partial quivers shown in Figure $7$.


\begin{figure}[htbp]

\beginpicture

\setcoordinatesystem units <0.8cm,0.6cm>             
\setplotarea x from -4 to 14, y from -2 to 6.5       

\put{Figure $7$: Partial quivers for $s_2s_3s_4s_3s_1s_2s_1s_3s_2s_4$.}[c] at 6 -3

\linethickness=0.5pt           

\put{$1$}[c] at -0.3 4 %
\put{$2$}[c] at -0.3 3 %
\put{$3$}[c] at -0.3 2 %
\put{$4$}[c] at -0.3 1 %
\put{$5$}[c] at -0.3 0 %

\put{$2$}[c] at 1.5 -1 %
\put{$3$}[c] at 2.5 -1 %
\put{$4$}[c] at 3.5 -1 %
\put{$3$}[c] at 4.5 -1 %
\put{$1$}[c] at 5.5 -1 %
\put{$2$}[c] at 6.5 -1 %
\put{$1$}[c] at 7.5 -1 %
\put{$3$}[c] at 8.5 -1 %
\put{$2$}[c] at 9.5 -1 %
\put{$4$}[c] at 10.5 -1 %

\setlinear \plot 0 0  3 0 / %
\setlinear \plot 3 0  5 2 / %
\setlinear \plot 5 2  6 2 / %
\setlinear \plot 6 2  8 4 / %
\setlinear \plot 8 4  12 4 / %

\setlinear \plot 0 1  2 1 / %
\setlinear \plot 2 1  3 2 / %
\setlinear \plot 3 2  4 2 / %
\setlinear \plot 4 2  5 1 / %
\setlinear \plot 5 1  8 1 / %
\setlinear \plot 8 1  10 3 / %
\setlinear \plot 10 3  12 3 / %

\setlinear \plot 0 2  1 2 / %
\setlinear \plot 1 2  2 3 / %
\setlinear \plot 2 3  5 3 / %
\setlinear \plot 5 3  6 4 / %
\setlinear \plot 6 4  7 4 / %
\setlinear \plot 7 4  8 3 / %
\setlinear \plot 8 3  9 3 / %
\setlinear \plot 9 3  10 2 / %
\setlinear \plot 10 2 12 2 / %

\setlinear \plot 0 3  1 3 / %
\setlinear \plot 1 3  4 0 / %
\setlinear \plot 4 0  10 0 / %
\setlinear \plot 10 0  11 1 / %
\setlinear \plot 11 1  12 1 / %

\setlinear \plot 0 4  5 4 / %
\setlinear \plot 5 4  7 2 / %
\setlinear \plot 7 2  8 2 / %
\setlinear \plot 8 2  9 1 / %
\setlinear \plot 9 1  10 1 / %
\setlinear \plot 10 1 11 0 / %
\setlinear \plot 11 0 12 0 / %

\put{$RRL$}[c] at 3.5 1.5 %
\put{$\0RL$}[c] at 4 2.5 %
\put{$\0\0L$}[c] at 7 0.5 %
\put{$LRL$}[c] at 6.5 1.5 %
\put{$\0R-$}[c] at 6.5 3.5 %
\put{$LR-$}[c] at 8 2.5 %

\endpicture

\end{figure}

We denote by $\mathcal{P}(\mathbf{i})$ the set of partial quivers obtained
from $\mathbf{i}$ in this way.

We now introduce the rectangles which we shall be considering. Let $i,j,k,l\in
\mathbb{N}$ satisfy:
$$i<j<l,\ \ i<k<l,\ \ i+l=j+k.$$
An $(i,j,k,l)$-rectangle is a rectangle with corners on levels $i,j,k,l$.
It is most convenient to illustrate this idea by means of an example. See
Figure $8$.


\begin{figure}[htbp]

\beginpicture

\setcoordinatesystem units <0.5cm,0.5cm>             
\setplotarea x from -11 to 9, y from 0 to 10.5      

\linethickness=1pt        

\put{Figure $8$: Drawing a $(0,2,3,5)$-rectangle.}[c] at 4 1.5

\setlinear \plot 3.5 9 	6.5 6 	/ %
\setlinear \plot 6.5 6 	4.5 4	/ %
\setlinear \plot 4.5 4 	1.5 7 	/ %
\setlinear \plot 1.5 7 	3.5 9 	/ %

\multiput {$0$}[c] at 1  9  *6 1 0 /      %
\multiput {$1$}[c] at 1  8  *6 1 0 /      %
\multiput {$2$}[c] at 1  7  *6 1 0 /      %
\multiput {$3$}[c] at 1  6  *6 1 0 /      %
\multiput {$4$}[c] at 1  5  *6 1 0 /      %
\multiput {$5$}[c] at 1  4  *6 1 0 /      %

\endpicture

\end{figure}

The sides of the rectangle have gradient $\pm \pi/4$. The
$(i,j,k,l)$-rectangle contains alternate columns of integers starting with
the first column if the entry in it is odd and the second column otherwise.
See Figure $9$.


\begin{figure}[htbp]

\beginpicture

\setcoordinatesystem units <0.5cm,0.5cm>             
\setplotarea x from -11 to 9, y from 0 to 10.5      

\linethickness=1pt        

\put{Figure $9$: A $(0,2,3,5)$-rectangle.}[c] at 5.25 1.5

\setlinear \plot 3.5 9 	6.5 6 	/ %
\setlinear \plot 6.5 6 	4.5 4	/ %
\setlinear \plot 4.5 4 	1.5 7 	/ %
\setlinear \plot 1.5 7 	3.5 9 	/ %

\put {$.$} at 2 7
\put {$1$} at 3 8
\put {$2$} at 3 7
\put {$3$} at 3 6
\put {$.$} at 4 8
\put {$.$} at 4 7
\put {$.$} at 4 6
\put {$.$} at 4 5
\put {$2$} at 5 7
\put {$3$} at 5 6
\put {$4$} at 5 5
\put {$.$} at 6 6

\endpicture

\end{figure}

The columns of integers in the rectangle are interpreted as positive roots;
for example the $(0,2,3,5)$-rectangle contains roots
$\alpha_1+\alpha_2+\alpha_3,\ \alpha_2+\alpha_3+\alpha_4$.

We next describe how each partial quiver determines a configuration of
rectangles. It is again most convenient to explain this by means of an
example. Consider the case of type $A_{10}$ in which we take the partial
quiver $P=\0LLRRRLRR$ --- see Figure $10$.

\begin{figure}[htbp]

\beginpicture

\setcoordinatesystem units <1cm,1cm>             
\setplotarea x from -1.8 to 12.5, y from 0 to 3       
\put {Figure $10$: The partial quiver $P=\0LLRRRLRR$.}[c] at 6.5 0.5

\scriptsize{

\multiput {$\circ$} at 2   2 *9 1 0 /      %

\linethickness=1pt           

\putrule from 2.05 2 to 2.95 2  %
\putrule from 3.05 2 to 3.95 2  %
\putrule from 4.05 2 to 4.95 2  %
\putrule from 5.05 2 to 5.95 2  %
\putrule from 6.05 2 to 6.95 2  %
\putrule from 7.05 2 to 7.95 2  %
\putrule from 8.05 2 to 8.95 2  %
\putrule from 9.05 2 to 9.95 2  %
\putrule from 10.05 2 to 10.95 2  %

\put {$1$}   [c] at 11.5 1.6
\put {$2$}   [c] at 10.5 1.6
\put {$3$}   [c] at 9.5 1.6
\put {$4$}   [c] at 8.5 1.6
\put {$5$}   [c] at 7.5 1.6
\put {$6$}   [c] at 6.5 1.6
\put {$7$}   [c] at 5.5 1.6
\put {$8$}   [c] at 4.5 1.6
\put {$9$}   [c] at 3.5 1.6
\put {$10$}  [c] at 2.5 1.6
\put {$11$}  [c] at 1.5 1.6

\setlinear \plot 3.6 2.1 3.5 2 / 
\setlinear \plot 4.6 2.1 4.5 2 / %
\setlinear \plot 8.6 2.1 8.5 2 / %

\setlinear \plot 3.6 1.9 3.5 2 / %
\setlinear \plot 4.6 1.9 4.5 2 / %
\setlinear \plot 8.6 1.9 8.5 2 / %

\setlinear \plot 5.4 2.1 5.5 2 / 
\setlinear \plot 6.4 2.1 6.5 2 / %
\setlinear \plot 7.4 2.1 7.5 2 / %
\setlinear \plot 9.4 2.1 9.5 2 / %
\setlinear \plot 10.4 2.1 10.5 2 / %

\setlinear \plot 5.4 1.9 5.5 2 / %
\setlinear \plot 6.4 1.9 6.5 2 / %
\setlinear \plot 7.4 1.9 7.5 2 / %
\setlinear \plot 9.4 1.9 9.5 2 / %
\setlinear \plot 10.4 1.9 10.5 2 / %

}

\endpicture

\end{figure}

The edges of the partial quiver are numbered as shown. We first divide the
partial quiver into its components, i.e. the maximal subquivers containing
a set of consecutive $L$'s and $R$'s. The components of our given partial
quiver $P$ are:
$$\begin{array}{ccccccccc}
 \0 &L & L & \0 & \0 & \0 & \0 & \0 & - \\
 \0 & \0 & \0 & R & R & R & \0 & \0 & -  \\
 \0 & \0 & \0 & \0 & \0 & \0 &L & \0 & - \\
 \0 & \0 & \0 & \0 & \0 & \0 & \0 &  R & R
\end{array}$$
For each component $K$ of $P$ we define positive integers $a(K)$, $b(K)$ with
$a(K)<b(K)$. The integer $a(K)$ is the number of the edge following the
rightmost arrow of
$K$ and $b(K)$ is the number of the edge preceding the leftmost arrow of $K$.
In the above example the numbers $a(K),b(K)$ are as follows:
$$\begin{array}{cccccccccccccc}
    &   &   &    &  K &    &    &    &   &   &   & a(K) & \ & b(K) \\
 \ & \ & \ & \ & \ & \ & \ & \ & \ & \ & \ & \ & \ & \ \\
 \0 & L & L & \0 & \0 & \0 & \0 & \0 & - & \ & \ & 7 & \ & 10 \\
 \0 & \0 & \0 & R & R & R & \0 & \0 & -  & \ & \ & 4 & \ & 8 \\
 \0 & \0 & \0 & \0 & \0 & \0 &L & \0 & - & \ & \ & 3 & \ & 5 \\
 \0 & \0 & \0 & \0 & \0 & \0 & \0 &  R & R & \ & \ & 1 & \ & 4
\end{array}$$
For each component $K$ of type $L$ we take a $(0,a,n+2-b,n+a-b+2)$-rectangle
and for each component $K$ of type $R$ we take a
$(b-a-1,b-1,n+1-a,n+1)$-rectangle, where $a=a(K)$, $b=b(K)$. Thus the $4$
components of our partial quiver $P=\0LLRRRLRR$ give the $4$ rectangles
shown in Figure $11$.

We now observe that for a component of type $L$ immediately followed by a
component of type $R$ the two corresponding rectangles fit together from
the left hand corner. Also, for a component of type $R$ immediately followed
by a component of type $L$ the two rectangles fit together from the right
hand corner. This can be observed in the four rectangles above, in which
the first and the second fit from the left, the second and the third from
the right, and the third and the fourth from the left.

We use these rules to superimpose the rectangles obtained from the components
of the given partial quiver. In the case of the partial quiver $P=
\0LLRRRLRR$ we obtain the configuration of rectangles shown in Figure $12$.


\begin{figure}[htbp]

\beginpicture

\setcoordinatesystem units <0.5cm,0.5cm>             
\setplotarea x from -3.5 to 10, y from -4 to 9      

\linethickness=1pt        

\put{$\0LL\0\0\0\0\0-$}[l] at 2 -1.7
\put{A $(0,7,2,9)$-rectangle}[l] at 2 -3

\setlinear \plot 0.5 2		7.5 9 / %
\setlinear \plot 7.5 9		9.5 7	 / %
\setlinear \plot 9.5 7		2.5 0  / %
\setlinear \plot 2.5 0	 	0.5 2 	 / %

\put {$7$} at 1 2
\multiput {$.$} at 2 1 *2 0 1 / %
\put {$8$} at 3 1
\put {$7$} at 3 2
\put {$6$} at 3 3
\put {$5$} at 3 4
\multiput {$.$} at 4 2 *3 0 1 / %
\put {$6$} at 5 3
\put {$5$} at 5 4
\put {$4$} at 5 5
\put {$3$} at 5 6
\multiput {$.$} at 6 4 *3 0 1 / %
\put {$4$} at 7 5
\put {$3$} at 7 6
\put {$2$} at 7 7
\put {$1$} at 7 8
\multiput {$.$} at 8 6 *2 0 1 / %
\put {$2$} at 9 7

\put{$\0\0\0RRR\0\0-$}[l] at 15 -1.7
\put{A $(3,7,7,11)$-rectangle}[l] at 15 -3

\setlinear \plot 14.5 4		18.5 8 / %
\setlinear \plot 18.5 8		22.5 4	 / %
\setlinear \plot 22.5 4		18.5 0  / %
\setlinear \plot 18.5 0	 	14.5 4 	 / %

\put {$7$} at 15 4
\multiput {$.$} at 16 3 *2 0 1 / %
\put {$9$} at 17 2
\put {$8$} at 17 3
\put {$7$} at 17 4
\put {$6$} at 17 5
\put {$5$} at 17 6
\multiput {$.$} at 18 1 *6 0 1 / %
\put {$10$} at 19 1
\put {$9$} at 19 2
\put {$8$} at 19 3
\put {$7$} at 19 4
\put {$6$} at 19 5
\put {$5$} at 19 6
\put {$4$} at 19 7
\multiput {$.$} at 20 2 *4 0 1 / %
\put {$8$} at 21 3
\put {$7$} at 21 4
\put {$6$} at 21 5
\put {$.$} at 22 4

\endpicture

\end{figure}


\begin{figure}[htbp]

\beginpicture

\setcoordinatesystem units <0.5cm,0.5cm>             
\setplotarea x from -3.5 to 10, y from -3.5 to 10      

\linethickness=1pt        

\put{$\0\0\0\0\0\0L\0-$}[l] at 2 -1.7
\put{A $(0,3,7,10)$-rectangle}[l] at 2 -3

\setlinear \plot 0.5 7		3.5 10 / %
\setlinear \plot 3.5 10		10.5 3	 / %
\setlinear \plot 10.5 3		7.5 0  / %
\setlinear \plot 7.5 0	 	0.5 7 	 / %

\put {$3$} at 1 7
\multiput {$.$} at 2 6 *2 0 1 / %
\put {$5$} at 3 5
\put {$4$} at 3 6
\put {$3$} at 3 7
\put {$2$} at 3 8
\put {$1$} at 3 9
\multiput {$.$} at 4 4 *5 0 1 / %
\put {$7$} at 5 3
\put {$6$} at 5 4
\put {$5$} at 5 5
\put {$4$} at 5 6
\put {$3$} at 5 7
\put {$2$} at 5 8
\multiput {$.$} at 6 2 *5 0 1 / %
\put {$9$} at 7 1
\put {$8$} at 7 2
\put {$7$} at 7 3
\put {$6$} at 7 4
\put {$5$} at 7 5
\put {$4$} at 7 6
\multiput {$.$} at 8 1 *4 0 1 / %
\put {$8$} at 9 2
\put {$7$} at 9 3
\put {$6$} at 9 4
\put {$.$} at 10 3

\put{$\0\0\0\0\0\0\0RR$}[l] at 15 -1.7
\put{A $(2,3,10,11)$-rectangle}[l] at 15 -3
\put{Figure $11$: The $4$ rectangles for $\0 LLRRRLRR$.}[c] at 13 -5.5

\setlinear \plot 14.5 8		15.5 9 / %
\setlinear \plot 15.5 9		23.5 1	 / %
\setlinear \plot 23.5 1		22.5 0  / %
\setlinear \plot 22.5 0	 	14.5 8 	 / %

\put {$3$} at 15 8
\multiput {$.$} at 16 7 *1 0 1 / %
\put {$5$} at 17 6
\put {$4$} at 17 7
\multiput {$.$} at 18 5 *1 0 1 / %
\put {$7$} at 19 4
\put {$6$} at 19 5
\multiput {$.$} at 20 3 *1 0 1 / %
\put {$9$} at 21 2
\put {$8$} at 21 3
\multiput {$.$} at 22 1 *1 0 1 / %
\put {$10$} at 23 1

\endpicture

\end{figure}


\begin{figure}[htbp]

\beginpicture

\setcoordinatesystem units <0.5cm,0.5cm>             
\setplotarea x from -11 to 10, y from -4 to 9.5      

\linethickness=1pt        

\put{Figure $12$: Configuration of all $4$ rectangles.}[c] at 5 -4


\setlinear \plot 0.5 2		7.5 9 / %
\setlinear \plot 7.5 9		9.5 7	 / %
\setlinear \plot 9.5 7		2.5 0  / %
\setlinear \plot 2.5 0	 	0.5 2 	 / %

\put {$7$} at 1 2
\multiput {$.$} at 2 1 *2 0 1 / %
\put {$8$} at 3 1
\put {$7$} at 3 2
\put {$6$} at 3 3
\put {$5$} at 3 4
\multiput {$.$} at 4 2 *3 0 1 / %
\put {$6$} at 5 3
\put {$5$} at 5 4
\put {$4$} at 5 5
\put {$3$} at 5 6
\multiput {$.$} at 6 4 *3 0 1 / %
\put {$4$} at 7 5
\put {$3$} at 7 6
\put {$2$} at 7 7
\put {$1$} at 7 8
\multiput {$.$} at 8 6 *2 0 1 / %
\put {$2$} at 9 7


\setlinear \plot 0.5 2		4.5 6 / %
\setlinear \plot 4.5 6		8.5 2	 / %
\setlinear \plot 8.5 2		4.5 -2  / %
\setlinear \plot 4.5 -2	 	0.5 2	 / %

\put {$7$} at 1 2
\multiput {$.$} at 2 1 *2 0 1 / %
\put {$9$} at 3 0
\put {$8$} at 3 1
\put {$7$} at 3 2
\put {$6$} at 3 3
\put {$5$} at 3 4
\multiput {$.$} at 4 -1 *6 0 1 / %
\put {$10$} at 5 -1
\put {$9$} at 5 0
\put {$8$} at 5 1
\put {$7$} at 5 2
\put {$6$} at 5 3
\put {$5$} at 5 4
\put {$4$} at 5 5
\multiput {$.$} at 6 0 *4 0 1 / %
\put {$8$} at 7 1
\put {$7$} at 7 2
\put {$6$} at 7 3
\put {$.$} at 7 2


\setlinear \plot -1.5 6		1.5 9 / %
\setlinear \plot 1.5 9		8.5 2	 / %
\setlinear \plot 8.5 2		5.5 -1  / %
\setlinear \plot 5.5 -1	 	-1.5 6 	 / %

\put {$3$} at -1 6
\multiput {$.$} at 0 5 *2 0 1 / %
\put {$5$} at 1 4
\put {$4$} at 1 5
\put {$3$} at 1 6
\put {$2$} at 1 7
\put {$1$} at 1 8
\multiput {$.$} at 2 3 *5 0 1 / %
\put {$7$} at 3 2
\put {$6$} at 3 3
\put {$5$} at 3 4
\put {$4$} at 3 5
\put {$3$} at 3 6
\put {$2$} at 3 7
\multiput {$.$} at 4 1 *5 0 1 / %
\put {$9$} at 5 0
\put {$8$} at 5 1
\put {$7$} at 5 2
\put {$6$} at 5 3
\put {$5$} at 5 4
\put {$4$} at 5 5
\multiput {$.$} at 6 0 *4 0 1 / %
\put {$8$} at 7 1
\put {$7$} at 7 2
\put {$6$} at 7 3
\put {$.$} at 8 2


\setlinear \plot -1.5 6		-0.5 7 / %
\setlinear \plot -0.5 7		7.5 -1	 / %
\setlinear \plot 7.5 -1		6.5 -2  / %
\setlinear \plot 6.5 -2		-1.5 6 	 / %

\put {$3$} at -1 6
\multiput {$.$} at 0 5 *1 0 1 / %
\put {$5$} at 1 4
\put {$4$} at 1 5
\multiput {$.$} at 2 3 *1 0 1 / %
\put {$7$} at 3 2
\put {$6$} at 3 3
\multiput {$.$} at 4 1 *1 0 1 / %
\put {$9$} at 5 0
\put {$8$} at 5 1
\multiput {$.$} at 6 -1 *1 0 1 / %
\put {$10$} at 7 -1

\endpicture

\end{figure}

We next define the {\em centre} of such a configuration. First consider the
rectangles in the diagonals from north-west to south-east. The number of such
rectangles in these diagonals in the above configuration is $1,3,4,2$. It is
always the case that one obtains a set of odd numbers followed by a set
of even numbers or vice versa. We draw the diagonal line separating the
diagonal blocks giving odd and even numbers of rectangles. This is the line
$\ell$ in Figure $13$. A similar phenomenon occurs for the diagonals from
northeast to southwest. The number of rectangles in such diagonals in the
above configuration is $2,4,3,1$ and we draw the diagonal line separating the
diagonal blocks giving odd and even numbers of rectangles. This is the line
$\ell'$ in Figure $13$. The point $O$ in which $\ell$ and $\ell'$
intersect is called the centre of the configuration and we draw in the
vertical line through $O$, called the central line, $m$. See Figure $13$.


\begin{figure}[htbp]

\beginpicture

\setcoordinatesystem units <0.5cm,0.5cm>             
\setplotarea x from -11.2 to 10, y from -5 to 10     

\linethickness=1pt        

\put{Figure $13$: Configuration of all $4$ rectangles with central line.}[c]
at 5 -5.3

\setlinear \plot -2 8.5  9 -2.5 / %
\setlinear \plot 0 -2.5  11 8.5 / %

\put {$\ell$} at -2.5 8
\put {$\ell$} at 8.5 -3
\put {$\ell'$} at 0.5 -3
\put {$\ell'$} at 11.5 8


\setlinear \plot 0.5 2		7.5 9 / %
\setlinear \plot 7.5 9		9.5 7	 / %
\setlinear \plot 9.5 7		2.5 0  / %
\setlinear \plot 2.5 0	 	0.5 2 	 / %

\put {$7$} at 1 2
\multiput {$.$} at 2 1 *2 0 1 / %
\put {$8$} at 3 1
\put {$7$} at 3 2
\put {$6$} at 3 3
\put {$5$} at 3 4
\multiput {$.$} at 4 2 *3 0 1 / %
\put {$6$} at 5 3
\put {$5$} at 5 4
\put {$4$} at 5 5
\put {$3$} at 5 6
\multiput {$.$} at 6 4 *3 0 1 / %
\put {$4$} at 7 5
\put {$3$} at 7 6
\put {$2$} at 7 7
\put {$1$} at 7 8
\multiput {$.$} at 8 6 *2 0 1 / %
\put {$2$} at 9 7


\setlinear \plot 0.5 2		4.5 6 / %
\setlinear \plot 4.5 6		8.5 2	 / %
\setlinear \plot 8.5 2		4.5 -2  / %
\setlinear \plot 4.5 -2	 	0.5 2	 / %

\put {$7$} at 1 2
\multiput {$.$} at 2 1 *2 0 1 / %
\put {$9$} at 3 0
\put {$8$} at 3 1
\put {$7$} at 3 2
\put {$6$} at 3 3
\put {$5$} at 3 4
\multiput {$.$} at 4 -1 *6 0 1 / %
\put {$10$} at 5 -1
\put {$9$} at 5 0
\put {$8$} at 5 1
\put {$7$} at 5 2
\put {$6$} at 5 3
\put {$5$} at 5 4
\put {$4$} at 5 5
\multiput {$.$} at 6 0 *4 0 1 / %
\put {$8$} at 7 1
\put {$7$} at 7 2
\put {$6$} at 7 3
\put {$.$} at 7 2


\setlinear \plot -1.5 6		1.5 9 / %
\setlinear \plot 1.5 9		8.5 2	 / %
\setlinear \plot 8.5 2		5.5 -1  / %
\setlinear \plot 5.5 -1	 	-1.5 6 	 / %

\put {$3$} at -1 6
\multiput {$.$} at 0 5 *2 0 1 / %
\put {$5$} at 1 4
\put {$4$} at 1 5
\put {$3$} at 1 6
\put {$2$} at 1 7
\put {$1$} at 1 8
\multiput {$.$} at 2 3 *5 0 1 / %
\put {$7$} at 3 2
\put {$6$} at 3 3
\put {$5$} at 3 4
\put {$4$} at 3 5
\put {$3$} at 3 6
\put {$2$} at 3 7
\multiput {$.$} at 4 1 *5 0 1 / %
\put {$9$} at 5 0
\put {$8$} at 5 1
\put {$7$} at 5 2
\put {$6$} at 5 3
\put {$5$} at 5 4
\put {$4$} at 5 5
\multiput {$.$} at 6 0 *4 0 1 / %
\put {$8$} at 7 1
\put {$7$} at 7 2
\put {$6$} at 7 3
\put {$.$} at 8 2


\setlinear \plot -1.5 6		-0.5 7 / %
\setlinear \plot -0.5 7		7.5 -1	 / %
\setlinear \plot 7.5 -1		6.5 -2  / %
\setlinear \plot 6.5 -2		-1.5 6 	 / %

\put {$3$} at -1 6
\multiput {$.$} at 0 5 *1 0 1 / %
\put {$5$} at 1 4
\put {$4$} at 1 5
\multiput {$.$} at 2 3 *1 0 1 / %
\put {$7$} at 3 2
\put {$6$} at 3 3
\multiput {$.$} at 4 1 *1 0 1 / %
\put {$9$} at 5 0
\put {$8$} at 5 1
\multiput {$.$} at 6 -1 *1 0 1 / %
\put {$10$} at 7 -1

\put {$m$} at 5 -3
\put {$m$} at 5 8
\setdashes <2mm,1mm>          %
\setlinear \plot 4.5 -3.5 4.5 8.5 / %


\put {$A$} at 10.2 7
\put {$B$} at 0 2
\put {$C$} at 9 2
\put {$D$} at -2 6
\put {$E$} at 8.2 -1

\endpicture

\end{figure}

We consider the left and right hand corner points in this configuration. These
are labelled $A,B,C,D,E$ in the above example. To each such corner point we
associate the rectangle which has the given point as a vertex and whose
edges through this point extend as far as possible in the figure.
(The vertex of the rectangle opposite to the given corner point may not be
explicitly shown in the figure). For each such corner point $V$ we define a
set of positive roots $\Phi^+(V)$. This is the set of all positive roots in
the rectangle associated with $V$ on the same side of the central line $m$ as
$V$ itself.

Thus in the given example, we have:
\begin{eqnarray*}
\Phi^+(A) & = & \{\alpha_2,\alpha_1+\alpha_2+\alpha_3+\alpha_4,
\alpha_3+\alpha_4+\alpha_5+\alpha_6\} \\
\Phi^+(B) & = & \{\alpha_7,\alpha_5+\alpha_6+\alpha_7+\alpha_8+\alpha_9\} \\
\Phi^+(C) & = & \{\alpha_6+\alpha_7+\alpha_8,
\alpha_4+\alpha_5+\alpha_6+\alpha_7+\alpha_8+\alpha_9+\alpha_{10}\} \\
\Phi^+(D) & = & \{\alpha_3,\alpha_1+\alpha_2+\alpha_3+\alpha_4+\alpha_5,
\alpha_2+\alpha_3+\alpha_4+\alpha_5+\alpha_6+\alpha_7\} \\
\Phi^+(E) & = & \{\alpha_{10},\alpha_8+\alpha_9\}
\end{eqnarray*}

We then define $\Phi^+(P)$ to be the union of the sets $\Phi^+(V)$ for all
corner points $V$ in the configuration. This is always a disjoint union. In
the given example we have:
\begin{eqnarray*}
\Phi^+(P) & = & \{
\alpha_2,\alpha_1+\alpha_2+\alpha_3+\alpha_4,
\alpha_3+\alpha_4+\alpha_5+\alpha_6,
\alpha_7,\alpha_5+\alpha_6+\alpha_7+\alpha_8+\alpha_9, \\
& & 
\alpha_6+\alpha_7+\alpha_8,
\alpha_4+\alpha_5+\alpha_6+\alpha_7+\alpha_8+\alpha_9+\alpha_{10},
\alpha_3,\alpha_1+\alpha_2+\alpha_3+\alpha_4+\alpha_5, \\
& & \alpha_2+\alpha_3+\alpha_4+\alpha_5+\alpha_6+\alpha_7,
\alpha_{10},\alpha_8+\alpha_9\}.
\end{eqnarray*}
We now define a vector $v_P\in \mathbb{N}^k$ whose coordinates are all $0$
or $1$. Let $\mathbf{j}$ be the reduced word
$$135\cdots 246\cdots 135\cdots 246\cdots, $$
considered above, and let $\alpha^1,\alpha^2,\ldots ,\alpha^k$ be the
corresponding order on the set of positive roots. We define $v_P$ as the
vector whose $i$th coordinate is $1$ if $\alpha^i\in \Phi^+(P)$ and is
$0$ otherwise. We also define vectors $v_j\in \mathbb{N}^k$ for $j=1,2,\ldots
,n$, where the $i$th component of $v_j$ is $1$ if the $i$th letter in
$\mathbf{j}$ is $j$ and is $0$ otherwise.

These vectors $v_P,\ P\in \mathcal{P}(\mathbf{i})$ and $v_j$, $j\in
\{1,2,\ldots ,n\}$ turn out to be our required spanning vectors for the
region $X_{\mathbf{i}}$.

\begin{prop}
Suppose that $\mathbf{g}$ has type $A_4$ and let $\mathbf{j}$ be the reduced
word $1324132413$ for $w_0$. Let $\mathbf{i}$ be any reduced word for $w_0$
and $\mathcal{P}(\mathbf{i})$ be the set of partial quivers associated with
$\mathbf{i}$. We have $|\mathcal{P}(\mathbf{i})|=6$. Then the region
$X_{\mathbf{i}}^+$ associated with $\mathbf{i}$ is the set of all non-negative
integral combinations of the vectors $v_P,\ P\in \mathcal{P}(\mathbf{i})$, and
$v_j$ for $1\leq j\leq 4$.
\end{prop}

This Proposition explains how the regions $X_{\mathbf{i}}^+$ can be
described by the rectangle combinatorics.

\section{Speculations for type $A_n$}

It is natural to ask whether the set $X_{\mathbf{i}}$ defined as the set of
non-negative linear combinations of the vectors $v_P$, $P\in \mathcal{P}
(\mathbf{i})$ and $v_j$, $1\leq j\leq n$ is a region of linearity for
$R:\mathbb{R}^k\rightarrow \mathbb{R}^k$ in type $A_n$. No counter-example
is known to the authors. If so, are the $X_{\mathbf{i}}$ the only regions of
linearity of $R$ defined by the minimum number of inequalities? This would
give a bijection in type $A_n$ between commutation classes of reduced words
for $w_0$ and regions of linearity for $R$ defined by the minimum number of
inequalities.

The set $M_{\mathbf{i}}$ of monomials corresponding to points in the Lusztig
cone $C_{\mathbf{i}}$ is not in general contained in the canonical basis in
type $A_n$. It is nevertheless possible to consider a subset of $\mathbf{B}$
in bijective correspondence with $C_{\mathbf{i}}$ by means of Kashiwara's
approach to the canonical basis. Kashiwara~\cite{kash2} defines certain
root operators $\widetilde{F_{\mathbf{i}}}$ which lead to a
parametrization of the canonical basis $\mathbf{B}$ for each $\mathbf{i}$ by
a certain subset $K_{\mathbf{i}}\subseteq \mathbb{N}^k$ which we call the
{\em string cone}. This gives a bijection
$$\psi_{\mathbf{i}}:\mathbf{B}\rightarrow K_{\mathbf{i}}.$$
It has been shown independently by Marsh~\cite{me9} and by
Premat~\cite{premat1}
that $C_{\mathbf{i}}\subseteq K_{\mathbf{i}}$ for each $\mathbf{i},$ i.e.
that the Lusztig cone lies in the string cone. Thus we obtain a subset
$\psi_{\mathbf{i}}^{-1}(C_{\mathbf{i}})\subseteq \mathbf{B}$. This subset
is equal to the set of monomials $M_{\mathbf{i}}$ when $n\leq 4$, but does
not consist of monomials in general. Using Lusztig's parametrization
$$\varphi_{\mathbf{j}}:\mathbf{B}\rightarrow \mathbb{N}^k$$
where $\mathbf{j}=135\cdots 246\cdots 135\cdots 246\cdots, $ this subset
$\psi_{\mathbf{i}}^{-1}(C_{\mathbf{i}})$ of $\mathbf{B}$ corresponds to a
certain subset of $\mathbb{N}^k$. Let
$$S_{\mathbf{i}}^{\mathbf{j}}:K_{\mathbf{i}}\rightarrow \mathbb{N}^k$$
be the transition map given by $S_{\mathbf{i}}^{\mathbf{j}}=\varphi_{\mathbf{j}}
\psi_{\mathbf{i}}^{-1}$. It is known that in type $A_n$, the images of the
spanning vectors of the Lusztig cone $C_{\mathbf{i}}$ under
$S_{\mathbf{i}}^{\mathbf{j}}$ are the vectors $v_P$,
$P\in \mathcal{P}(\mathbf{i})$ and $v_j$, $1\leq j\leq n$ given by the
rectangle combinatorics. This can be proved using a transition function
introduced by Berenstein, Fomin and Zelevinsky~\cite{bfz1}.

Finally, what can be said about canonical basis elements corresponding
under $\varphi_{\mathbf{j}}$ to regions of linearity of $R$ not defined by
the minimum number of inequalities? The results of Xi~\cite{xi3} in type
$A_3$ are interesting in this respect. In type $A_3$, Lusztig's function
$R:\mathbb{R}^6\rightarrow \mathbb{R}^6$ has $10$ regions of linearity,
$8$ of which are the regions $X_{\mathbf{i}}$ for the different commutation
classes of reduced words $\mathbf{i}$ for $w_0$. These are all defined
by $3$ inequalities. The remaining two are defined by $4$ inequalities and
we denote these by $X_9$ and $X_{10}$. We now define
$\widetilde{X_{\mathbf{i}}^+}$ to be the set of points in $X_{\mathbf{i}}$
whose coordinates are all real and nonnegative. Thus we have
$$X_{\mathbf{i}}^+\subseteq \widetilde{X_{\mathbf{i}}^+}\subseteq
X_{\mathbf{i}}$$
and $X_{\mathbf{i}}^+$ is the set of integral points in
$\widetilde{X_{\mathbf{i}}^+}$. We define
$\widetilde{X_9^+}$ and $\widetilde{X_{10}^+}$ similarly. We have
additional inequalities defining
$\widetilde{X_{\mathbf{i}}^+}$ asserting that all coordinates are
non-negative, but some of these inequalities will be redundant. In fact each
of the $8$ regions $\widetilde{X_{\mathbf{i}}^+}$ can be defined by $6$
inequalities and, for suitable numbering, $\widetilde{X_9^+}$ can be defined
by $8$ inequalities and $\widetilde{X_{10}^+}$ by $9$ inequalities. The
regions $\widetilde{X_{\mathbf{i}}^+}$ are called simplicial regions as the
number of defining inequalities is equal to the dimension of the ambient
space. $\widetilde{X_9^+}$ can be written as the union of two simplicial
regions, and $\widetilde{X_{10}^+}$ as the union of four simplicial regions.

Xi obtains $8$ families of monomials in $\mathbf{B}$, which are parametrized
by the integral points in the $8$ simplicial regions
$\widetilde{X_{\mathbf{i}}^+}$. In addition he obtains $6$ families of
elements in $\mathbf{B}$ which are not monomials. These are parametrized by
the integral points in the two simplicial regions whose union is
$\widetilde{X_9^+}$ and the four simplicial regions whose union is
$\widetilde{X_{10}^+}$.
The canonical basis elements corresponding to a given region are all of the
same type, i.e. they can all be expressed as a linear combination of
monomials corresponding to a fixed reduced expression with quantum binomial
coefficients and exponents lying on a line segment in $\mathbb{N}^6$.

The authors plan to give the proofs of the results described in this article
in a forthcoming paper.~\nocite{me10}

\newcommand{\noopsort}[1]{}\newcommand{\singleletter}[1]{#1}

\end{document}